\documentclass[12pt]{article}

\usepackage{amsmath}
\usepackage{amssymb}
\usepackage{amscd}
\usepackage{amsthm}
\usepackage{indentfirst}
\usepackage[english,frenchb]{babel}

\newtheorem{thm}{Th\'eor\`eme}[section]
\newtheorem{lem}[thm]{Lemme}

\newtheorem{prop}[thm]{Proposition}

\theoremstyle{remark}
\newtheorem{rmq}[thm]{Remarque}
\newtheorem{exm}[thm]{Exemple}

\DeclareMathOperator{\Cl}{Cl}

\DeclareMathOperator{\Spec}{Spec}

\DeclareMathOperator{\pic}{Pic}
\DeclareMathOperator{\ord}{ord}

\DeclareMathOperator{\rpic}{RPic}

\newcommand{\gm}{\mathbf{G}_{{\rm m}}}

\newcommand{\gmss}{\mathbf{G}_{{\rm m},S'}}

\newcommand{\znz}{(\mathbb{Z}/n\mathbb{Z})_S}
\DeclareMathOperator{\im}{Im}
\newcommand{\Hom}{\underline{\rm Hom}}

\newcommand{\homr}{{\rm Hom}}
\newcommand{\Ext}{\underline{\rm Ext}^1}
\newcommand{\ext}{{\rm Ext}^1}
\newcommand{\extp}{{\rm Ext}_P}

\newcommand{\res}{\mathrm{res}_K}
\newcommand{\He}{\mathcal{H}}

\begin{document}

\title{Invariants de classes : propri\'et\'es fonctorielles et applications \`a l'\'etude du noyau}

\author{Jean Gillibert}

\date{septembre 2007}

\maketitle

\begin{abstract}
L'homomorphisme de classes \'etudie la structure galoisienne de torseurs -- sous un sch\'ema en groupes fini et plat -- obtenus gr\^ace au cobord d'une suite exacte. Son introduction est due \`a Martin Taylor (la suite exacte \'etant une isog\'enie entre sch\'emas ab\'eliens). Nous commen\c cons par \'enoncer quelques propri\'et\'es g\'en\'erales de cet homomorphisme, puis nous poursuivons son \'etude dans le cas o\`u la suite exacte est donn\'ee par la multiplication par $n$ sur une extension d'un sch\'ema ab\'elien par un tore.

\medskip

\begin{otherlanguage}{english}
\begin{center}
\textbf{Abstract}
\end{center}

The class-invariant homomorphism allows one to measure the Galois module structure of torsors---under a finite flat group scheme---which lie in the image of a coboundary map associated to an exact sequence. It has been introduced first by Martin Taylor (the exact sequence being given by an isogeny between abelian schemes). We begin by giving general properties of this homomorphism, then we pursue its study in the case when the exact sequence is given by the multiplication by $n$ on an extension of an abelian scheme by a torus.
\end{otherlanguage}
\end{abstract}



\section{Introduction}

Soit $F/K$ une extension galoisienne de corps de nombres, et soit $\Gamma$ son groupe de Galois.
L'action de $\Gamma$ sur $F$ donne naissance \`a une structure de $K[\Gamma]$-module sur $F$. Le th\'eor\`eme de la base normale affirme que $F$ est libre de rang $1$ en tant $K[\Gamma]$-module.

Soit \`a pr\'esent $\mathcal{O}_K$ (resp. $\mathcal{O}_F$) l'anneau des entiers de $K$ (resp. $F$). On constate que $\Gamma$ agit encore sur $\mathcal{O}_F$, ce qui munit ce dernier d'une structure de $\mathcal{O}_K[\Gamma]$-module. L'\'etude de cette structure a \'et\'e initi\'ee par D. Hilbert dans son {\it Zahlbericht} \cite{ref16}.

Le crit\`ere de Noether affirme que $\mathcal{O}_F$ est un
$\mathcal{O}_K[\Gamma]$-module projectif (de type fini) si et seulement
si l'extension $F/K$ est mod\'er\'ement ramifi\'ee. Dans ce cas, on peut \'etudier la classe de $\mathcal{O}_F$ dans le groupe des classes localement libres $\Cl(\mathcal{O}_K[\Gamma])$ (ce dernier, que l'on note \'egalement $\tilde{\mathrm{K}}_0(\mathcal{O}_K[\Gamma])$, est par d\'efinition le noyau du morphisme $\mathrm{K}_0(\mathcal{O}_K[\Gamma])\rightarrow \mathbb{Z}$ qui \`a la classe d'un module associe son rang local). La conjecture de Fr\"ohlich, montr\'ee par Taylor \cite{ref20}, est le point culminant de cette th\'eorie.

Par contre, lorsque l'extension est sauvagement ramifi\'ee, la situation est moins claire (voir \cite[Appendix, C]{ref7}). Une approche alternative au probl\`eme initial consiste \`a remplacer $\mathcal{O}_K[\Gamma]$ par une certaine alg\`ebre de Hopf.

On suppose \`a pr\'esent que $\Gamma$ est commutatif. Alors l'ensemble des classes d'isomorphie de
$K$-alg\`ebres galoisiennes de groupe $\Gamma$ forme un groupe, qui n'est autre que le groupe de cohomologie (galoisienne, \'etale ou fppf, selon les pr\'ef\'erences du lecteur) $H^1(K,\Gamma)$.

Soit $\He$ un $\mathcal{O}_K$-ordre de Hopf dans $K[\Gamma]$, c'est-\`a-dire que $\He$ est une $\mathcal{O}_K$-alg\`ebre de Hopf, projective de type fini en tant que $\mathcal{O}_K$-module, telle que $\He\otimes_{\mathcal{O}_K} K\simeq K[\Gamma]$. Comme pr\'ec\'edemment, l'ensemble des classes d'isomorphie de $\He$-extensions galoisiennes forme un groupe, qui s'identifie au groupe de cohomologie fppf $H^1(\Spec(\mathcal{O}_K),\Spec(\He^*))$, o\`u $\He^*$ est l'alg\`ebre de Hopf duale de $\He$. En outre, la fl\`eche naturelle de restriction \`a la fibre g\'en\'erique
$$
\begin{CD}
\res:H^1(\Spec(\mathcal{O}_K),\Spec(\He^*)) @>>> H^1(K,\Gamma) \\
\end{CD}
$$
est injective. Ainsi, si l'on choisit $\He$ de telle sorte que la classe de $F$ soit dans l'image de $\res$, on obtient un unique $\mathcal{O}_K$-ordre $\mathcal{N}$ dans $F$ qui est muni d'une structure de $\He$-extension galoisienne. En particulier, $\mathcal{N}$ est muni d'une action de $\He$. On constate que $\mathcal{N}$ est $\He$-localement libre de rang $1$, donc d\'efinit une classe $(\mathcal{N})$ dans le groupe $\Cl(\He)$. On montre alors que l'application
\begin{equation*}
\begin{split}
\pi_0:H^1(\Spec(\mathcal{O}_K),\Spec(\He^*))\;\longrightarrow & \;\Cl(\He)\\
\mathcal{N}\;\longmapsto & \;(\mathcal{N})(\He^*)^{-1}\\
\end{split}
\end{equation*}
est un morphisme de groupes.
On dit que $\mathcal{N}$ admet une base normale si $\pi_0(\mathcal{N})=0$, ce qui revient \`a dire que $\mathcal{N}$ et $\He^*$ sont isomorphes en tant que $\He$-modules.

\subsection{Structure galoisienne des torseurs}

Soit $S$ un sch\'ema n\oe th\'erien, et soit $G$ un $S$-sch\'ema en groupes commutatif, fini et plat. Nous noterons $G^D$ le dual de Cartier de $G$.

Waterhouse a d\'efini dans \cite[Theorem $5$]{ref22} un morphisme de groupes
$$
\begin{CD}
\pi:H^1(S,G) @>>> \pic(G^D)\\
\end{CD}
$$
qui g\'en\'eralise l'homomorphisme $\pi_0$ d\'efini ci-dessus, le groupe de classes $\Cl(\He)$ pouvant \^etre identifi\'e au groupe de Picard $\pic(\Spec(\He))$ (ceci est d\^u au fait que la dimension de Krull de $\mathcal{O}_K$, donc de $\He$, est \'egale \`a $1$, voir \cite[Cor. 3.5]{ref4}). Nous dirons que $\pi$ mesure la structure galoisienne des $G$-torseurs. Ainsi un $G$-torseur a une structure galoisienne triviale, ou, de fa\c con \'equivalente, admet une base normale, si son image par $\pi$ est nulle.

Une fa\c con conceptuelle de d\'efinir $\pi$ est la suivante. On d\'eduit de la suite spectrale \og locale-globale \fg{} pour les Ext \cite[expos\'e V, proposition 6.3, 3)]{ref14} la suite exacte
\begin{align*}
0\rightarrow H^1(S,\Hom(G^D,\gm)) \rightarrow \ext(G^D,\gm) \rightarrow \Ext(G^D,\gm)(S) \rightarrow \cdots
\end{align*}
o\`u $\gm$ d\'esigne le groupe multiplicatif sur $S$. Le faisceau $\Hom_S(G^D,\gm)$ est repr\'esentable par $G$, et le faisceau $\Ext(G^D,\gm)$ est nul d'apr\`es \cite[expos\'e VIII, 3.3.1]{ref15}. Par suite, on dispose d'un isomorphisme
\begin{equation}
\label{chelou}
\begin{CD}
H^1(S,G) @>\sim>> \ext(G^D,\gm).\\
\end{CD}
\end{equation}
On d\'efinit alors $\pi$ comme \'etant \'egal au compos\'e de cet isomorphisme avec le morphisme naturel $\ext(G^D,\gm)\rightarrow \pic(G^D)$.

\subsection{Probl\`eme central}
\label{pbcentral}

Consid\'erons une suite exacte de la forme
\begin{equation}
\label{mainsuite}
\begin{CD}
0 @>>> N(f) @>>> F_1 @>f>> F_2 @>>> 0\\
\end{CD}
\end{equation}
o\`u $F_1$ et $F_2$ sont deux faisceaux ab\'eliens pour la topologie fppf sur $S$, tels que le faisceau noyau $N(f)$ soit repr\'esent\'e par un $S$-sch\'ema en groupes fini et plat. Gr\^ace au foncteur des sections globales, on en d\'eduit un morphisme cobord
$$
\begin{CD}
\delta:F_2(S)@>>> H^1(S,N(f))
\end{CD}
$$
et l'on souhaite \'etudier la structure galoisienne des torseurs ainsi obtenus.

Plus formellement, on d\'efinit un homomorphisme $\psi_f$ comme \'etant obtenu par composition du cobord $\delta$ et de l'homomorphisme $\pi$
$$
\begin{CD}
\psi_f:F_2(S) @>\delta>> H^1(S,N(f)) @>\pi>> \pic(N(f)^D).\\
\end{CD}
$$
On dit, pour reprendre l'appellation usuelle, que $\psi_f$ est l'homomorphisme de classes associ\'e \`a la suite exacte (\ref{mainsuite}).

\subsection{Pr\'esentation des r\'esultats}

Dans la section 2 de cet article, nous \'enon\c cons quelques propri\'et\'es fonctorielles (comportement par changement de base, par produit, etc.) de $\pi$ et de $\psi_f$. Ces propri\'et\'es, bien qu'\'el\'emen\-taires, ne semblent pas \^etre apparues pour l'instant dans la litt\'erature, c'est pourquoi nous les avons incluses ici.

Ces pr\'eliminaires nous permettent, dans la section 3, d'\'etudier le noyau de l'homomorphisme de classes associ\'e \`a la suite exacte de Kummer (multiplication par un entier $n>0$) dans une extension d'un sch\'ema ab\'elien par un tore.

Les notations introduites ici seront en vigueur tout au long de cet article.
Nous fixons, comme pr\'ec\'edemment, un sch\'ema n\oe th\'erien $S$. Nous notons $\gm$ le groupe multiplicatif sur $S$. De plus, pour tout entier naturel $n$, nous notons $\mu_n$ le sch\'ema en groupes des racines $n$-i\`emes de l'unit\'e sur $S$.

Dans toute la suite, nous adoptons la convention suivante : un rev\^etement est un morphisme fini localement libre de rang constant. Si $S'\rightarrow S$ est un rev\^etement, nous noterons $[S':S]$ son rang.

Par d\'efinition, un $S$-tore est un $S$-sch\'ema en groupes qui est localement isomorphe, pour la topologie \'etale sur $S$, \`a une puissance de $\gm$. On dit qu'un tore est d\'eploy\'e s'il est isomorphe \`a une puissance de $\gm$.

Si $X$ est une extension d'un $S$-sch\'ema ab\'elien par un $S$-tore, alors pour tout entier $n>0$ la multiplication par $n$ dans $X$ est un \'epimorphisme fppf, et son noyau est un $S$-sch\'ema en groupes fini et plat. Nous noterons $\psi_n^X$ l'homomorphisme de classes associ\'e.

Nous commen\c cons par \'etudier l'homomorphisme de classes associ\'e \`a un $S$-tore. Nous montrons dans le paragraphe \ref{kertor} le r\'esultat suivant.

\begin{thm}
\label{nultor}
Soit $T$ un $S$-tore.
\begin{enumerate}
\item[$(i)$] Si $T$ est d\'eploy\'e, alors $\psi_n^T$ est nul pour tout entier $n$.
\item[$(ii)$] Supposons que $S$ soit connexe normal, et soit $S'\rightarrow S$ un rev\^etement (\'etale) sur lequel $T$ se d\'eploye. Alors $\psi_n^T$ est nul pour tout entier $n$ premier \`a $[S':S]$.
\end{enumerate}
\end{thm}

Clarifions la situation consid\'er\'ee dans $(ii)$ : si $S$ est normal, alors $T$ est isotrivial d'apr\`es \cite[expos\'e X, th\'eor\`eme 5.16]{ref13}, c'est-\`a-dire qu'il existe un morphisme \'etale fini $S'\rightarrow S$ qui d\'eploye $T$. Si en outre $S$ est connexe, alors $S'\rightarrow S$ est de rang constant, et par cons\'equent est un rev\^etement dans notre terminologie.

\begin{rmq}
Dans le paragraphe \ref{cextor}, nous expliquons comment construire des suites exactes de la forme $(\ref{mainsuite})$ o\`u $F_1$ et $F_2$ sont des tores, et telles que $\psi_f$ ne soit pas nul. L'hypoth\`ese sur $n$ dans le cas $(ii)$ ci-dessus semble donc \^etre n\'ecessaire.
\end{rmq}

Puis nous \'etudions l'homomorphisme de classes associ\'e \`a un $S$-sch\'ema ab\'elien. On sait que, si $E$ est une $S$-courbe elliptique et si $n$ est premier \`a $6$, alors $\psi_n^E$ s'annule sur les points de torsion. Ce r\'esultat, conjectur\'e par Taylor dans \cite{ref21}, a \'et\'e montr\'e en premier par Taylor et Srivastav dans \cite{ref19}, puis g\'en\'eralis\'e par Agboola dans \cite{ref2}, et enfin par Pappas dans \cite{ref18}. Pour une interpr\'etation arithm\'etique (resp. g\'eom\'etrique) de ce probl\`eme, on pourra consulter \cite{ref6} (resp. \cite{ref1}). Dans le paragraphe \ref{kerab}, nous \'etendons ce r\'esultat de la fa\c con suivante.

\begin{thm}
\label{nulisoell}
Soit $A$ un $S$-sch\'ema ab\'elien. Supposons que l'une des conditions suivantes soit v\'erifi\'ee :
\begin{enumerate}
\item[$(i)$] Il existe une isog\'enie $\nu:A\rightarrow C$ de degr\'e $N$, telle que $C$ soit un produit de $S$-courbes elliptiques.
\item[$(ii)$] Il existe un rev\^etement $S'\rightarrow S$ de rang $N$ tel que $A_{S'}$ soit isomorphe \`a un produit de $S'$-courbes elliptiques.
\end{enumerate}
Alors les points de torsion de $A(S)$ sont dans le noyau de $\psi_n^A$ pour tout entier $n$ premier \`a $6N$.
\end{thm}

\begin{rmq}
Signalons ici que, dans le cas de figure $(ii)$, la conclusion du th\'eor\`eme n'est plus vraie si l'on supprime l'hypoth\`ese que $n$ est premier \`a $6N$. Plus pr\'ecis\'ement, on exhibe dans \cite[exemple 4.10]{ref10} un sch\'ema ab\'elien $A$ de dimension $5$, muni d'un point de $5$-torsion qui n'appartient pas au noyau de  $\psi_5^A$. De plus, $A$ est obtenu comme restriction de Weil d'une courbe elliptique, et en tant que tel devient isomorphe \`a un produit de courbes elliptiques apr\`es changement de base par un rev\^etement de rang $5$.
\end{rmq}

Enfin, dans le paragraphe \ref{dernier}, nous consid\'erons une extension d'un $S$-sch\'ema ab\'elien par le groupe multiplicatif.

\begin{thm}
\label{nulpourext}
Soient $A$ un $S$-sch\'ema ab\'elien, et $V$ une extension de $A$ par $\gm$, dont nous noterons $f:V\rightarrow A$ la projection naturelle. Supposons que l'une des conditions suivantes soit v\'erifi\'ee :
\begin{enumerate}
\item[$(i)$] L'entier $n$ est premier \`a $|\pic(S)|$ et $V[n]^D(S)$ est de cardinal \'egal \`a son ordre.
\item[$(ii)$] La classe de $V$ dans le groupe $\ext(A,\gm)$ est un \'el\'ement d'ordre fini premier \`a $n$.
\end{enumerate}
Alors $f^{-1}(\ker(\psi_n^A))=\ker(\psi_n^V)$.
\end{thm}

La deuxi\`eme assertion de ce th\'eor\`eme se r\'ev\`ele particuli\`erement utile pour construire de nouveaux exemples d'annulation de l'homomorphisme de classes \`a partir d'une courbe elliptique, comme nous le montrons dans l'exemple qui suit.

\begin{exm}
Soit $E$ une $S$-courbe elliptique, et soit $t\in E(S)$ un point de torsion d'ordre $N$. Alors $t$ d\'etermine, en vertu de l'isomorphisme d'auto-dualit\'e
$E(S)\simeq \ext(E,\gm)$, une extension $V_t$ de $E$ par $\gm$, dont nous noterons $f_t:V_t\rightarrow E$ la projection naturelle. Soit $n$ un entier premier \`a $6N$, on consid\`ere l'homomorphisme de classes
$$
\begin{CD}
\psi_n:V_t(S) @>>> \pic(V_t[n]^D).\\
\end{CD}
$$
En combinant le th\'eor\`eme \ref{pappastaylor} et le $(ii)$ du th\'eor\`eme \ref{nulpourext}, on obtient le r\'esultat suivant :
$$f_t^{-1}(E(S)_{\rm Tors})\subseteq \ker(\psi_n).$$
Le cas le plus int\'eressant est quand on part d'un point de torsion $x\in E(S)$ d'ordre premier \`a $N$. Alors l'image r\'eciproque de $x$ par $f_t$ est non vide (plus exactement est en bijection avec $\Gamma(S,\mathcal{O}_S)^{\times}$), et fournit une nouvelle famille de torseurs dont la structure galoisienne est triviale.
\end{exm}


\section{Propri\'et\'es fonctorielles}

Dans toute cette section, $G$ d\'esigne un $S$-sch\'ema en groupes commutatif, fini et plat, et $G^D$ d\'esigne son dual de Cartier. Nous noterons $\ord(G)$ l'ordre de $G$, c'est-\`a-dire le rang du morphisme fini localement libre $G\rightarrow S$.

\subsection{Comportement par changement de base}
\label{weilrestriction}

Soit  $h:S'\rightarrow S$ un rev\^etement de rang $[S':S]$.
Si $H'$ est un $S'$-sch\'ema en groupes affine, le faisceau $h_*H'$ est repr\'esentable par un $S$-sch\'ema en groupes \'egalement affine, que l'on appelle la restriction de Weil de $H'$, et que l'on note $\Re_{S'/S}(H')$. Pour les d\'etails, voir \cite[\S\/7.6]{ref5}.

Soit \`a pr\'esent $H$ un $S$-sch\'ema en groupes affine commutatif. Soit
$$N_{S'/S}:\Re_{S'/S}(H_{S'})\rightarrow H$$
le morphisme \og trace \fg{} \cite[expos\'e XVII, 6.3.13]{ref14}. Rappelons que la compos\'ee
$$
\begin{CD}
H @>>> \Re_{S'/S}(H_{S'}) @>N_{S'/S}>> H \\
\end{CD}
$$
est la multiplication par $[S':S]$ dans le groupe $H$ \cite[expos\'e XVII, 6.3.15]{ref14}. Supposons que $H$ soit lisse sur $S$, on dispose alors d'un isomorphisme
\begin{equation}
\label{johnb}
H^1(S,\Re_{S'/S}(H_{S'}))\simeq H^1(S',H_{S'}).
\end{equation}
En effet, les groupes $H_{S'}$ et $\Re_{S'/S}(H_{S'})$ \'etant lisses, les deux $H^1$ sont les m\^emes pour la topologie \'etale ou la topologie fppf. On se sert alors du fait que $R^1h_*(H_{S'})=0$ en topologie \'etale, le morphisme $h$ \'etant fini (voir \cite[expos\'e VIII, 5.5 et 5.3]{ref14}).

Dans le cas o\`u $H=\gm$, le morphisme \og trace \fg{} $\Re_{S'/S}(\gmss)\rightarrow \gm$ donne ainsi naissance \`a un morphisme
$$\mathbf{N}_{S'/S}:\pic(S')\longrightarrow \pic(S)$$
que nous appellerons morphisme \og norme \fg{}. Par abus de notation, pour tout $S$-sch\'ema $X$ nous d\'esignerons \'egalement par $\mathbf{N}_{S'/S}:\pic(X_{S'})\longrightarrow \pic(X)$ le morphisme norme associ\'e au rev\^etement $X_{S'}\rightarrow X$. Pour une autre d\'efinition du morphisme norme, on pourra consulter \cite[\S\/6.5]{ref12}.

Par fonctorialit\'e du morphisme $\pi$, on constate la commutativit\'e du diagramme
\begin{equation}
\label{pi2}
\begin{CD}
H^1(S',G_{S'}) @>\pi'>> \pic(G_{S'}^D)\\
@AAA @AAA \\
H^1(S,G) @>\pi>> \pic(G^D)\\
\end{CD}
\end{equation}
o\`u les fl\`eches verticales sont obtenues par changement de base.

Signalons le r\'esultat suivant, dont nous ne nous servirons pas dans la suite, mais qui pr\'esente un int\'er\^et en lui-m\^eme.

\begin{lem}
Supposons que $[S':S]$ soit premier \`a l'ordre de $G$. Alors la fl\`eche de changement de base
$$H^1(S,G) \longrightarrow H^1(S',G_{S'})$$
est injective. En d'autres termes, aucun $G$-torseur n'est trivialis\'e par un rev\^etement de rang premier \`a l'ordre de $G$.
\end{lem}

\begin{proof}
Nous avons une suite spectrale
$$H^p(S,R^qh_*(G_{S'}))\Longrightarrow H^{p+q}(S',G_{S'})$$
dont on d\'eduit une suite exacte
$$
0\longrightarrow H^1(S,\Re_{S'/S}(G_{S'})) \longrightarrow H^1(S',G_{S'}) \longrightarrow H^0(S,R^1h_*(G_{S'})).\\
$$
D'autre part, nous avons une fl\`eche $H^1(S,G)\rightarrow H^1(S,\Re_{S'/S}(G_{S'}))$ induite par l'application canonique $G\rightarrow \Re_{S'/S}(G_{S'})$. Cette fl\`eche est injective. En effet, la compos\'ee
$$
\begin{CD}
H^1(S,G) @>>> H^1(S,\Re_{S'/S}(G_{S'})) @>>> H^1(S,G)\\
\end{CD}
$$
est \'egale \`a la multiplication par $[S':S]$ dans $H^1(S,G)$, laquelle est un automorphisme de $H^1(S,G)$, puisque $[S':S]$ est premier \`a l'ordre de $G$.

Au final, le morphisme de changement de base $H^1(S,G) \rightarrow H^1(S',G_{S'})$, qui s'\'ecrit comme la compos\'ee des fl\`eches
$$
\begin{CD}
H^1(S,G) @>>> H^1(S,\Re_{S'/S}(G_{S'})) @>>> H^1(S',G_{S'}) \\
\end{CD}
$$
est injectif. Ce qu'on voulait.
\end{proof}

\begin{prop}
\label{chb}
Soient $Z\in H^1(S,G)$ un $G$-torseur, et $Z_{S'}\in H^1(S',G_{S'})$ le torseur obtenu par changement de base.
\begin{enumerate}
\item[$(i)$] La relation $\pi'(Z_{S'})=0$ implique $\pi(Z)^{[S':S]}=0$.
\item[$(ii)$] Supposons que $[S':S]$ soit premier \`a l'ordre de $G$. Alors $\pi(Z)$ est nul si et seulement si $\pi'(Z_{S'})$ l'est.
\end{enumerate}
\end{prop}

\begin{proof}
$(i)$ Consid\'erons le morphisme \og norme \fg{}
$$\mathbf{N}_{S'/S}:\pic(G_{S'}^D)\longrightarrow \pic(G^D)$$
et rappelons que, pour tout $\mathcal{L}\in \pic(G^D)$, nous avons $\mathbf{N}_{S'/S}(\mathcal{L}_{S'})=\mathcal{L}^{[S':S]}$.
Supposons \`a pr\'esent que $Z$ satisfasse $\pi'(Z_{S'})=0$. Alors, par commutativit\'e du diagramme (\ref{pi2}), nous avons $\pi(Z)_{S'}=0$. Par suite,
$\mathbf{N}_{S'/S}(\pi(Z)_{S'})=\pi(Z)^{[S':S]}=0$, ce qu'on voulait.

$(ii)$ La condition est n\'ecessaire par commutativit\'e du diagramme (\ref{pi2}). Etablissons \`a pr\'esent la suffisance. Le groupe $H^1(S,G)$ \'etant tu\'e par l'entier $\ord(G)$, nous avons $\pi(Z)^{\ord(G)}=0$. De plus, on sait que $\pi(Z)^{[S':S]}=0$ d'apr\`es le point $(i)$. Les entiers $[S':S]$ et $\ord(G)$ \'etant premiers entre eux, il en r\'esulte que $\pi(Z)=0$, ce qu'on voulait.
\end{proof}

Soit \`a pr\'esent une suite exacte de la forme
$$
\begin{CD}
0 @>>> G @>>> F_1 @>f>> F_2 @>>> 0\\
\end{CD}
$$
dans laquelle $F_1$ et $F_2$ sont des faisceaux ab\'eliens pour la topologie fppf sur $S$. Comme convenu dans le paragraphe \ref{pbcentral}, nous notons $\psi_f:F_2(S)\rightarrow \pic(G^D)$ l'homomorphisme de classes associ\'e \`a cette suite.

D'autre part, le foncteur $h^*$ est exact (consid\'er\'e en tant que foncteur de la cat\'egorie des faisceaux ab\'eliens pour la topologie fppf sur $S$ dans la cat\'egorie analogue sur $S'$). On d\'eduit donc de la suite pr\'ec\'edente une suite exacte
$$
\begin{CD}
0 @>>> G_{S'} @>>> h^*F_1 @>h^*f>> h^*F_2 @>>> 0\\
\end{CD}
$$
et l'on note $\psi_{h^*f}:F_2(S')\rightarrow \pic(G_{S'}^D)$ l'homomorphisme associ\'e \`a cette suite.
\begin{prop}
\label{chbpsi}
Soit $x\in F_2(S)$, et soit $x'\in F_2(S')$ la $S'$-section associ\'ee.
\begin{enumerate}
\item[$(i)$] La relation $\psi_{h^*f}(x')=0$ implique $\psi_f(x)^{[S':S]}=0$.
\item[$(ii)$] Supposons que $[S':S]$ soit premier \`a l'ordre de $G$. Alors $\psi_f(x)$ est nul si et seulement si $\psi_{h^*f}(x')$ l'est.
\end{enumerate}
\end{prop}

\begin{proof}
Le diagramme suivant est commutatif
$$
\begin{CD}
F_2(S') @>\delta'>> H^1(S',G_{S'}) \\
@AAA @AAA \\
F_2(S) @>\delta >> H^1(S,G) \\
\end{CD}
$$
o\`u $\delta$ (resp. $\delta'$) d\'esigne le cobord associ\'e \`a la suite exacte de d\'epart (resp. \`a la suite exacte obtenue par changement de base), et o\`u les fl\`eches verticales sont obtenues par changement de base. La proposition \ref{chb} permet d'en d\'eduire le r\'esultat.
\end{proof}


\subsection{Comportement fonctoriel}

Consid\'erons un diagramme commutatif \`a lignes exactes de la forme
$$
\begin{CD}
0 @>>> N(f) @>>> F_1 @>f>> F_2 @>>> 0 \\
@. @. @V\nu_1VV @V\nu_2VV \\
0 @>>> N(g) @>>> H_1 @>g>> H_2 @>>> 0 \\
\end{CD}
$$
o\`u $F_1$, $F_2$, $H_1$ et $H_2$ sont des faisceaux ab\'eliens pour la topologie fppf sur $S$. Alors la fl\`eche $\nu_1$ induit une fl\`eche $\nu_0:N(f)\rightarrow N(g)$. Supposons que $N(f)$ et $N(g)$ soient repr\'esent\'es par des sch\'emas en groupes finis et plats sur $S$. On peut alors \'enoncer la proposition suivante.

\begin{prop}
\label{foncteur}
Soit $\nu_0^D:N(g)^D\rightarrow N(f)^D$ la fl\`eche d\'eduite de $\nu_0$ par dualit\'e de Cartier. Alors le diagramme
$$
\begin{CD}
F_2(S) @>\psi_f>> \pic(N(f)^D) \\
@V\nu_2VV @VV(\nu_0^D)^*V \\
H_2(S) @>\psi_g>> \pic(N(g)^D) \\
\end{CD}
$$
est commutatif.
\end{prop}

\begin{proof}
Tout d'abord, le diagramme suivant (obtenu par passage \`a la cohomologie dans les deux suites exactes de d\'epart)
$$
\begin{CD}
@>>> F_1(S) @>f>> F_2(S) @>\delta_f>> H^1(S,N(f)) @>>> \cdots \\
@. @V\nu_1VV @V\nu_2VV @VV(\nu_0)_*V \\
@>>> H_1(S) @>g>> H_2(S) @>\delta_g>> H^1(S,N(g)) @>>> \cdots \\
\end{CD}
$$
est commutatif. De plus, le diagramme
$$
\begin{CD}
H^1(S,N(f))\, @. \simeq \ext(N(f)^D,\gm) @>\pi_f>> \pic(N(f)^D) \\
@VV(\nu_0)_*V @VVV @VV(\nu_0^D)^*V \\
H^1(S,N(g))\, @. \simeq \ext(N(g)^D,\gm) @>\pi_g>> \pic(N(g)^D)  \\
\end{CD}
$$
est \'egalement commutatif. D'o\`u le r\'esultat.
\end{proof}

\begin{prop}
\label{foncteurcor}
Supposons que $\nu_0:N(f)\rightarrow N(g)$ soit un isomorphisme. Alors
$$
\ker(\psi_f)=\nu_2^{-1}(\ker(\psi_g)).
$$
\end{prop}

\begin{proof}
Soit $x\in F_2(S)$. Nous allons montrer que $\psi_g(\nu_2(x))$ est nul si et seulement si $\psi_f(x)$ l'est. D'apr\`es la proposition \ref{foncteur}, la condition $\psi_g(\nu_2(x))=0$ \'equivaut \`a la condition $\nu_0^*(\psi_f(x))=0$. Le morphisme $\nu_0$ \'etant un isomorphisme, on en d\'eduit que $\nu_0^*$ en est \'egalement un, donc la condition pr\'ec\'edente se r\'ecrit $\psi_f(x)=0$, ce qu'on voulait.
\end{proof}


\subsection{Comportement par produit}

Consid\'erons deux suites exactes
$$
\begin{CD}
0 @>>> N(f) @>>> F_1 @>f>> F_2 @>>> 0\\
\end{CD}
$$
et
$$
\begin{CD}
0 @>>> N(g) @>>> H_1 @>g>> H_2 @>>> 0\\
\end{CD}
$$
o\`u $F_1$, $F_2$, $H_1$ et $H_2$ sont des faisceaux ab\'eliens pour la topologie fppf sur $S$. On suppose que $N(f)$ et $N(g)$ sont repr\'esent\'es par des sch\'emas en groupes finis et plats sur $S$.

Soit $\psi_f$ (resp. $\psi_g$) l'homomorphisme associ\'e \`a la premi\`ere (resp. deux\-i\`eme) suite, et soit $\psi_{f\times g}$ l'homomorphisme associ\'e \`a la suite produit
$$
\begin{CD}
0 @>>> N(f)\times_S N(g) @>>> F_1\times_S H_1 @>f\times g>> F_2\times_S H_2 @>>> 0.\\
\end{CD}
$$
Nous avons alors la propri\'et\'e suivante :

\begin{prop}
\label{produit}
Avec les notations pr\'ec\'edentes, on a l'\'egalit\'e
$$
\ker \psi_{f\times g}=\ker \psi_f\times \ker \psi_g.
$$
\end{prop}

\begin{proof}
Consid\'erons le diagramme (commutatif) suivant
$$
\begin{CD}
H^1(S,N(f))\times H^1(S,N(g)) @>\pi_f\times\pi_g>> \pic(N(f)^D)\times \pic(N(g)^D)\\
@VVV @VVV \\
H^1(S,N(f)\times_S N(g)) @>\pi_{f\times g}>> \pic(N(f)^D\times_S N(g)^D)\\
\end{CD}
$$
dans lequel la fl\`eche de gauche est un isomorphisme, que l'on obtient en composant les isomorphismes qui suivent
\begin{equation*}
\begin{split}
H^1(S,N(f))\times H^1(S,N(g)) & \simeq \ext(N(f)^D,\gm)\times
\ext( N(g)^D,\gm) \\
& \simeq \ext(N(f)^D\times_S N(g)^D,\gm) \\
& \simeq H^1(S,N(f)\times_S N(g)).\\
\end{split}
\end{equation*}
Il suffit de montrer que, sous cet isomorphisme, $\ker \pi_f\times \ker \pi_g$ s'identifie \`a $\ker \pi_{f\times g}$. L'inclusion $\ker \pi_f\times \ker \pi_g\subseteq\ker \pi_{f\times g}$ est imm\'ediate au vu du diagramme qui pr\'ec\`ede. D'autre part, la fl\`eche de droite dans le susdit diagramme se d\'ecrit explicitement de la fa\c con suivante
\begin{equation*}
\begin{split}
\pic(N(f)^D)\times \pic(N(g)^D)\;\longrightarrow & \;\pic(N(f)^D\times_S N(g)^D)\\
(\mathcal{L}_1,\mathcal{L}_2)\;\longmapsto & \;pr_1^*(\mathcal{L}_1)+pr_2^*(\mathcal{L}_2)\\
\end{split}
\end{equation*}
o\`u $pr_1$ et $pr_2$ sont les projections canoniques du produit $N(f)^D\times_S N(g)^D$ sur ses deux facteurs. Mais $N(f)^D$ et $N(g)^D$ sont des $S$-sch\'emas en groupes ; gr\^ace \`a leurs sections neutres on obtient respectivement des sections $s_1$ et $s_2$ de $pr_1$ et $pr_2$. Il est alors clair que l'application $(s_1^*,s_2^*)$ constitue une section du morphisme ci-dessus, lequel est donc injectif. On en d\'eduit ais\'ement l'inclusion $\ker \pi_{f\times g}\subseteq\ker \pi_f\times \ker \pi_g$, d'o\`u le r\'esultat.
\end{proof}


\subsection{Comportement sur certaines suites exactes}

Si $G$ est un $S$-sch\'ema en groupes fini et plat, nous notons $\rpic(G^D)$ l'image de l'homomorphisme de Waterhouse ; on dispose d'une suite exacte
\begin{equation}
\label{triw}
\begin{CD}
0 \rightarrow \extp^1(G^D,\gm) @>>> H^1(S,G) @>\pi >> \rpic(G^D) \rightarrow 0\\
\end{CD}
\end{equation}
o\`u $\extp^1(G^D,\gm)$ est le groupe des extensions de $G^D$ par $\gm$ dans la cat\'egorie des pr\'efaisceaux sur $S$. On appelle usuellement $\rpic(G^D)$ le groupe des classes r\'ealisables. Nous renvoyons \`a Agboola \cite{ref3} pour de plus amples d\'etails sur cette suite exacte.

On consid\`ere \`a pr\'esent une suite exacte de $S$-sch\'emas en groupes fini et plats pour la topologie fppf sur $S$
\begin{equation}
\label{base}
\begin{CD}
0 @>>> G_0 @>>> G_1 @>>> G_2 @>>> 0\\
\end{CD}
\end{equation}
et l'on note respectivement $\pi_0$, $\pi_1$ et $\pi_2$ les homomorphismes de Waterhouse associ\'es aux groupes $G_0$, $G_1$ et $G_2$.

\begin{lem}
\label{pfxsuite}
Supposons que la suite (obtenue \`a partir de (\ref{base}) par dualit\'e de Cartier)
\begin{equation}
\label{dbase}
\begin{CD}
0 @>>> G_2^D @>>> G_1^D @>>> G_0^D @>>> 0\\
\end{CD}
\end{equation}
soit une suite exacte dans la cat\'egorie des pr\'efaisceaux sur $S$. Alors on dispose d'une suite exacte
$$
\begin{CD}
0 @>>> \rpic(G_0^D) @>>> \rpic(G_1^D) @>>> \rpic(G_2^D)\\
\end{CD}
$$
dans la cat\'egorie des groupes ab\'eliens.
\end{lem}

\begin{proof}
La suite (\ref{triw}) est fonctorielle en $G$, de sorte que nous pouvons d\'eduire de la suite exacte (\ref{base}) un diagramme commutatif :
\[
\begin{CD}
0\longrightarrow\ @.\homr(G_2^D,\gm) @>\sim>> G_2(S) @>>> 0\\
@. @VVV @VVV @VVV \\
0\longrightarrow\ @.\extp^1(G_0^D,\gm) @>>> H^1(S,G_0) @>\pi_0>> \rpic(G_0^D)@.\ \longrightarrow 0 \\
@. @VVV @VVV @VVV \\
0\longrightarrow\ @.\extp^1(G_1^D,\gm) @>>> H^1(S,G_1) @>\pi_1>> \rpic(G_1^D)@.\ \longrightarrow 0 \\
@. @VVV @VVV @VVV \\
0\longrightarrow\ @.\extp^1(G_2^D,\gm) @>>> H^1(S,G_2) @>\pi_2>> \rpic(G_2^D)@.\ \longrightarrow 0 \\
\end{CD}
\]
dans lequel les trois derni\`eres lignes sont des suites exactes d\'eduites de (\ref{triw}). La colonne du milieu est exacte, par passage \`a la cohomologie dans la suite (\ref{base}). De plus, l'exactitude de la suite (\ref{dbase}) dans la cat\'egorie des pr\'efaisceaux implique l'exactitude de la colonne de gauche.
Nous modifions \`a pr\'esent ce diagramme en supprimant la premi\`ere ligne et en rempla\c cant respectivement dans la deuxi\`eme ligne les deux premiers termes par les conoyaux des fl\`eches $\homr(G_2^D,\gm)\rightarrow\extp^1(G_0^D,\gm)$ et $G_2(S)\rightarrow H^1(S,G_0)$. Dans le diagramme ainsi obtenu, les lignes, ainsi que les deux premi\`eres colonnes, sont exactes. D'apr\`es le lemme des neuf appliqu\'e \`a ce nouveau diagramme, la troisi\`eme colonne, qui n'a pas \'et\'e alt\'er\'ee par les modifications effectu\'ees, est une suite exacte. Ce qu'on voulait.
\end{proof}

\begin{rmq}
Si les facteurs de composition des fibres de $G$ au-dessus des points de $S$ de caract\'eristique r\'esiduelle $2$ ne contiennent pas de facteur $\alpha_2$, alors d'apr\`es \cite[Theorem 1.2]{ref3} nous avons
$$
\rpic(G^D)=\homr(G^D,H^1(-,\gm)),
$$
\'egalit\'e dans laquelle $H^1(-,\gm)$ est consid\'er\'e en tant que foncteur en groupes ab\'eliens, ou pr\'efaisceau. Si l'on suppose que $G_0$, $G_1$ et $G_2$ satisfont aux susdites hypoth\`eses, alors notre lemme d\'ecoule de ce r\'esultat en appliquant le foncteur $\homr(-,H^1(-,\gm))$ \`a la suite exacte (\ref{dbase}) dans la cat\'egorie des pr\'efaisceaux.

D'autre part, le pr\'efaisceau $H^1(-,\gm)$ n'est pas un faisceau pour la topologie fppf. Il n'y a donc aucune raison pour que $\rpic$ soit  exact \`a gauche en g\'en\'eral.
\end{rmq}


\section{Applications \`a l'\'etude du noyau}

\subsection{Tores : groupe multiplicatif et cons\'equences}
\label{kertor}

Soit $n>0$ un entier naturel fix\'e. Nous avons une suite exacte de faisceaux ab\'eliens (pour la topologie fppf sur $S$)
\[
\begin{CD}
0 @>>> \mu_n @>>> \gm @>[n]>> \gm @>>> 0.\\
\end{CD}
\]
Soit $d$ le cobord associ\'e. On obtient, par application du foncteur des sections globales, une suite exacte
\[
\begin{CD}
\cdots \longrightarrow \gm(S) @>d>> H^1(S,\mu_n) @>s>> H^1(S,\gm)[n] \longrightarrow 0.\\
\end{CD}
\]
On montre ais\'ement le r\'esultat suivant :

\begin{prop}
\label{nulgm}
Soit $\pi$ le morphisme de Waterhouse pour le groupe $\mu_n$. Alors, avec les notations pr\'ec\'edentes, $\im d=\ker \pi$.
\end{prop}

\begin{proof}
On peut donner une description explicite de $H^1(S,\mu_n)$ (voir \cite[page 125]{ref17}). Plus pr\'ecis\'ement, $H^1(S,\mu_n)$ s'identifie \`a l'ensemble des (classes d'isomorphie de) couples $(\mathcal{L}, \sigma)$, o\`u $\mathcal{L}$ est un $\gm$-torseur sur $S$, et $\sigma: \mathcal{L}^{\otimes n}\rightarrow\mathcal{O}_S$ est un isomorphisme.
Avec cette description, $d$ est l'application qui \`a un \'el\'ement $\alpha\in \gm(S)=\Gamma(S,\mathcal{O}_S)^{\times}$ associe le couple $(\mathcal{O}_S,\sigma)$, o\`u $\sigma:\mathcal{O}_S^{\otimes n}\simeq \mathcal{O}_S\rightarrow \mathcal{O}_S$ est la multiplication par $\alpha$. De plus, $s$ est l'application qui au couple $(\mathcal{L}, \sigma)$ associe $\mathcal{L}\in\pic(S)$.

Soit $\mathcal{M}(\mathcal{L},\sigma)$ la $\mathcal{O}_S$-alg\`ebre dont le spectre est le $\mu_n$-torseur correspondant au couple $(\mathcal{L}, \sigma)$. En tant que $\mathcal{O}_S$-module, on peut d\'ecrire $\mathcal{M}(\mathcal{L},\sigma)$ de la fa\c con suivante
$$\mathcal{M}(\mathcal{L},\sigma)=\oplus_{k=0}^{n-1}\mathcal{L}^{\otimes k}$$
o\`u l'on a pos\'e $\mathcal{L}^{\otimes 0}=\mathcal{O}_S$ par convention. La loi de multiplication de $\mathcal{M}(\mathcal{L},\sigma)$ est induite par le produit tensoriel et par l'application $\sigma$.

Soit \`a pr\'esent $C_n$ un groupe (abstrait) cyclique d'ordre $n$, de g\'en\'erateur $g$ et d'\'el\'ement neutre $e$. Alors $\mu_n$ est le spectre de l'alg\`ebre de groupe $\mathcal{O}_S[C_n]$. L'action de $\mu_n$ sur le torseur $(\mathcal{L}, \sigma)$ correspond \`a une coaction de $\mathcal{O}_S[C_n]$ sur $\mathcal{M}(\mathcal{L},\sigma)$, que l'on peut d\'ecrire explicitement de la fa\c con suivante :
\begin{equation*}
\begin{split}
\oplus_{k=0}^{n-1}\mathcal{L}^{\otimes k} \;\longrightarrow & \;
\mathcal{O}_S[C_n]\otimes_{\mathcal{O}_S} (\oplus_{k=0}^{n-1}\mathcal{L}^{\otimes k})\\
(s_0,s_1,\dots,s_{n-1}) \;\longmapsto & \; (e\otimes s_0, g\otimes s_1,\dots, g^{n-1}\otimes s_{n-1}).\\
\end{split}
\end{equation*}
Soit $\mathcal{O}_S[C_n]^*$ l'alg\`ebre duale de $\mathcal{O}_S[C_n]$, dont le spectre est le $S$-sch\'ema en groupes constant $(\mathbb{Z}/n\mathbb{Z})_S$. Alors l'action de $\mathcal{O}_S[C_n]^*$ sur $\mathcal{M}(\mathcal{L},\sigma)$ induite par la coaction qui pr\'ec\`ede est tout simplement la multiplication composante par composante. En d'autres termes, pour tout $f\in \mathcal{O}_S[C_n]^*$ et pour tout $(s_0,s_1,\dots,s_{n-1})\in \mathcal{M}(\mathcal{L},\sigma)$,
$$f\cdot(s_0,s_1,\dots,s_{n-1})=(f(e)s_0,f(g)s_1,\dots,f(g^{n-1})s_{n-1}).$$

On peut en d\'eduire que le morphisme de Waterhouse
$$\pi:H^1(S,\mu_n)\longrightarrow \pic((\mathbb{Z}/n\mathbb{Z})_S)\simeq \pic(S)^n$$
est l'application qui \`a $(\mathcal{L}, \sigma)$ associe le $n$-uplet $(\mathcal{O}_S,\mathcal{L},\mathcal{L}^{\otimes 2},\dots, \mathcal{L}^{\otimes (n-1)})$ dans le groupe $\pic(S)^n$.
Il est donc clair que $\pi((\mathcal{L}, \sigma))=0$ si et seulement si $\mathcal{L}$ est nul dans $\pic(S)$, d'o\`u le r\'esultat.
\end{proof}

\begin{rmq}
\label{immediate}
Les morphismes $\pi$ et $s$ ayant m\^eme noyau (\`a savoir $\im d$), leurs images $\rpic(\znz)$ et $\pic(S)[n]$ sont isomorphes.
\end{rmq}

\begin{proof}[D\'emonstration du th\'eor\`eme \ref{nultor}]
$(i)$ La proposition \ref{nulgm} implique que l'homomorphisme de classes associ\'e \`a la multiplication par $n$ dans $\gm$ est nul. Le m\^eme r\'esultat est donc \'egalement valable pour n'importe quelle puissance de $\gm$, en vertu de la proposition \ref{produit}. $(ii)$ Il suffit de combiner le $(i)$ et la proposition \ref{chbpsi}.
\end{proof}

\begin{rmq}
Supposons que $n$ soit inversible sur $S$, et que $\mathcal{O}_S(S)$ contienne une racine primitive $n$-i\`eme de l'unit\'e. Alors $\mu_n$ est isomorphe, en tant que $S$-sch\'ema en groupes, au sch\'ema en groupes constant $(\mathbb{Z}/n\mathbb{Z})_S$. Ainsi la proposition 6.5 de \cite[chap. 0]{ref11} d\'ecoule de notre proposition \ref{nulgm}.
\end{rmq}


\subsection{Tores : contre-exemples en dimension sup\'erieure}
\label{cextor}

Nous donnons \`a pr\'esent un proc\'ed\'e de construction de suites exactes de la forme $(\ref{mainsuite})$ o\`u $F_1$ et $F_2$ sont des tores, et telles que $\psi_f$ soit non nul. Ce processus a \'et\'e adapt\'e au cas des vari\'et\'es ab\'eliennes dans \cite{ref10}.

Soient $K$ un corps de nombres et $\mathcal{O}_K$ l'anneau des entiers de $K$. Dans ce paragraphe, notre sch\'ema de base sera $S=\Spec(\mathcal{O}_K)$.

\begin{prop}
\label{beautore}
Soit $n>0$ un entier. Alors il existe un $S$-tore $T$ et une suite exacte
\begin{equation}
\begin{CD}
0 @>>> \mu_n @>>> T @>f>> Q @>>> 0\\
\end{CD}
\end{equation}
tels que le cobord qui s'en d\'eduit
$$
\begin{CD}
\delta: Q(S) @>>> H^1(S,\mu_n) \\
\end{CD}
$$
soit surjectif. Par suite, l'image de l'homomorphisme de classes $\psi_f$ correspondant est le groupe $\rpic(\znz)=\pic(S)[n]$ tout entier.
\end{prop}

\begin{rmq}
Ce r\'esultat implique que l'homomorphisme de classes associ\'e \`a une isog\'enie entre tores n'est pas nul en g\'en\'eral. Il suffit en effet de choisir $K$ et $n$ de telle mani\`ere que $n$ divise l'ordre de $\pic(S)$.
\end{rmq}

\begin{proof}
Soit $K'$ une extension finie de $K$, non ramifi\'ee en toutes les places finies de $K$, et soit $S'=\Spec(\mathcal{O}_{K'})$ le spectre de l'anneau des entiers de $K'$. On suppose en outre que l'application naturelle induite par le changement de base $h:S'\rightarrow S$
$$
\begin{CD}
h^*_n:\pic(S)[n] @>>> \pic(S')[n] \\
\end{CD}
$$
est nulle. Une telle extension $K'$ existe d'apr\`es la th\'eorie du corps de classes (tout id\'eal de $K$ se principalise dans le corps de classes de Hilbert de $K$, qui est une extension non ramifi\'ee de $K$).

Soit \`a pr\'esent le sch\'ema $T:=\Re_{S'/S}(\gmss)$. Alors $T$ est un $S$-tore car $T$ est trivialis\'e par le rev\^etement \'etale $S'\rightarrow S$. D'autre part, nous avons un morphisme canonique
$$
\begin{CD}
c_1:\gm @>>> T\\
\end{CD}
$$
qui est une immersion ferm\'ee d'apr\`es \cite[\S\/7.6, p. 197]{ref5}, le groupe $\gm$ \'etant affine, donc s\'epar\'e sur $S$. Ainsi, $\mu_n$ est un sous-sch\'ema en groupes de $T$. Soit $Q$ le quotient (pour la topologie fppf) de $T$ par $\mu_n$ ($Q$ est \'egalement repr\'esentable par un $S$-tore), de sorte que nous avons une suite exacte
$$
\begin{CD}
0 @>>> \mu_n @>>> T @>f>> Q @>>> 0.\\
\end{CD}
$$

Consid\'erons \`a pr\'esent le diagramme commutatif
$$
\begin{CD}
0 @>>> \mu_n @>>> \gm @>[n]>> \gm @>>> 0\\
@. @| @Vc_1VV @Vc_2VV \\
0 @>>> \mu_n @>>> T @>f>> Q @>>> 0\\
\end{CD}
$$
dans lequel l'application $c_2$ est d\'eduite de $c_1$ par passage au quotient. Par passage \`a la cohomologie, et en se servant du fait que le groupe $H^1(S,\mu_n)$ est tu\'e par $n$, nous obtenons un diagramme commutatif \`a lignes exactes
$$
\begin{CD}
\gm(S) @>d>> H^1(S,\mu_n) @>>> H^1(S,\gm)[n] \\
@Vc_2VV @| @VVV \\
Q(S) @>\delta>> H^1(S,\mu_n) @>>> H^1(S,T)[n]. \\
\end{CD}
$$
D'autre part, le groupe $\gm$ \'etant lisse sur $S$, nous avons, gr\^ace \`a l'isomorphisme (\ref{johnb}) du paragraphe \ref{weilrestriction}, une suite d'identifications
$$
H^1(S,T)=H^1(S,\Re_{S'/S}(\gmss))\simeq H^1(S',\gmss)=\pic(S'),
$$
de sorte que le diagramme pr\'ec\'edent se traduit par le diagramme suivant (\`a lignes exactes, comme son pr\'ed\'ecesseur)
$$
\begin{CD}
\gm(S) @>d>> H^1(S,\mu_n) @>>> \pic(S)[n] \\
@Vc_2VV @| @VVh^*_nV \\
Q(S) @>\delta>> H^1(S,\mu_n) @>q>> \pic(S')[n] \\
\end{CD}
$$
dans lequel la fl\`eche de droite n'est autre que l'application $h^*_n$ d\'eduite du changement de base $S'\rightarrow S$. Or, par hypoth\`ese, cette application est nulle. On en d\'eduit aussit\^ot, par commutativit\'e du carr\'e droit du diagramme, que la fl\`eche $q$ est nulle, donc que le morphisme $\delta$ est surjectif, ce qu'on voulait.
\end{proof}


\subsection{Vari\'et\'es ab\'eliennes : courbes elliptiques et cons\'equences}
\label{kerab}

Nous renvoyons le lecteur \`a l'article de Pappas \cite{ref18} pour une d\'emonstration du r\'esultat suivant.

\begin{thm}
\label{pappastaylor}
Supposons que $E$ soit une $S$-courbe elliptique et que $n$ soit un entier premier \`a $6$. Alors les points de torsion de $E(S)$ sont dans le noyau de $\psi_n^E$.
\end{thm}

Comme nous l'avons soulign\'e dans l'introduction, le th\'eor\`eme \ref{nulisoell} est une cons\'equence imm\'ediate de ce r\'esultat.

\begin{proof}[D\'emonstration du th\'eor\`eme \ref{nulisoell}]
$(i)$ Soit $\psi_n^C$ l'homomorphisme de classes associ\'e \`a la multiplication par $n$ dans $C$. Sachant que $C$ est un produit de courbes elliptiques et que $n$ est premier \`a $6$, on d\'eduit ais\'ement du th\'eor\`eme \ref{pappastaylor}, avec l'aide de la proposition \ref{produit}, que $\psi_n^C$ s'annule sur les points de torsion de $C(S)$, autrement dit $C(S)_{\rm Tors}\subseteq \ker(\psi_n^C)$. Comme l'image par $\nu$ d'un point de torsion est un point de torsion, on en d\'eduit que
$$
\nu(A(S)_{\rm Tors})\subseteq \ker(\psi_n^C).
$$
D'autre part, $n$ \'etant premier \`a $N$, l'homomorphisme $\nu_0:A[n]\rightarrow C[n]$ obtenu par restriction de $\nu$ est un isomorphisme. D'apr\`es la proposition \ref{foncteurcor}, nous avons donc
$$
\ker(\psi_n^A)=\nu^{-1}(\ker(\psi_n^C))
$$
ce qui, combin\'e \`a l'inclusion pr\'ec\'edente, permet de conclure.
$(ii)$ Il suffit d'appliquer la proposition \ref{chbpsi}, combin\'ee au th\'eor\`eme \ref{pappastaylor}.
\end{proof}


\subsection{Extensions de vari\'et\'es ab\'eliennes par le groupe multiplicatif}
\label{dernier}

Soient $A$ un $S$-sch\'ema ab\'elien et $n>0$ un entier naturel. Consid\'erons une extension
$$
\begin{CD}
0 @>>> \gm @>>> V @>f>> A @>>> 0\\
\end{CD}
$$
de $A$ par $\gm$. Par des arguments standards on en d\'eduit une suite exacte
\begin{equation}
\label{tors}
\begin{CD}
0 @>>> \mu_n @>>> V[n] @>>> A[n] @>>> 0.\\
\end{CD}
\end{equation}

Nous noterons $\psi_n^V$ (resp. $\psi_n^A$) l'homomorphisme de classes associ\'e \`a la multiplication par $n$ dans $V$ (resp. $A$). En vertu de la proposition \ref{foncteur}, nous obtenons un diagramme commutatif
\begin{equation}
\label{foncteurplus}
\begin{CD}
V(S) @>\psi_n^V >> \rpic(V[n]^D) \\
@VfVV @VVV \\
A(S) @>\psi_n^A >> \rpic(A[n]^D). \\
\end{CD}
\end{equation}
Une question naturelle se pose alors : quelles relations peut-on donner entre le noyau de $\psi_n^V$ et celui de $\psi_n^A$ ? Le th\'eor\`eme \ref{nulpourext} constitue un premier r\'esultat dans cette direction.

\begin{lem}
\label{clef}
On suppose que $V[n]^D(S)$ est de cardinal \'egal \`a l'ordre de $V[n]^D$. Alors la suite (d\'eduite de $(\ref{tors})$ par dualit\'e de Cartier)
\begin{equation}
\label{dualtors}
\begin{CD}
0 @>>> A[n]^D @>>> V[n]^D @>>> \znz @>>> 0\\
\end{CD}
\end{equation}
est exacte dans la cat\'egorie des pr\'efaisceaux ab\'eliens sur $S$.
\end{lem}

\begin{proof}
Pour montrer cela, il faut montrer que $V[n]^D$ est trivial en tant que $A[n]^D$-torseur sur $\znz$. Comme $\znz$ est isomorphe (en tant que sch\'ema) \`a la r\'eunion disjointe de $n$ copies de $S$, $V[n]^D$ est \'egal \`a la r\'eunion disjointe de $n$ $A[n]^D$-torseurs sur $S$. De plus, chacun de ces torseurs admet une section sur $S$ (car $V[n]^D$ admet autant de sections sur $S$ que son ordre). Donc chacun de ces torseurs est trivial, d'o\`u le r\'esultat.
\end{proof}

\begin{proof}[D\'emonstration du th\'eor\`eme \ref{nulpourext}]
$(i)$ D'apr\`es le lemme \ref{clef}, la suite (\ref{dualtors}) est exacte dans la cat\'egorie des pr\'efaisceaux ab\'eliens sur $S$. D'apr\`es le lemme \ref{pfxsuite} on obtient une suite exacte
$$
\begin{CD}
0 \longrightarrow \rpic(\znz) @>>> \rpic(V[n]^D) @>>> \rpic(A[n]^D).\\
\end{CD}
$$
Or on sait que le premier terme de cette suite est isomorphe \`a $\pic(S)[n]$, donc est nul d'apr\`es les hypoth\`eses. Autrement dit, la fl\`eche de droite dans le diagramme (\ref{foncteurplus}) est injective. Le r\'esultat en d\'ecoule ais\'ement.

$(ii)$ Soit $N$ l'ordre de la classe de $V$ dans $\ext(A,\gm)$. On sait que la multiplication par $N$ dans le groupe $\ext(A,\gm)$ est induite par la multiplication par $N$ dans $\gm$. Donc, si l'on note $NV$ l'extension repr\'esent\'ee par la suite du bas,
$$
\begin{CD}
0 @>>> \gm @>>> V @>f>> A @>>> 0\\
@. @V[N]VV @VV\phi V @| \\
0 @>>> \gm @>>> NV @>g>> A @>>> 0\\
\end{CD}
$$
laquelle est obtenue par push-out de la suite du haut, alors la classe de $NV$ est isomorphe \`a la classe triviale, en particulier $NV\simeq \gm\times_S A$ en tant que faisceaux. Soit $\psi_n^{NV}$ l'homomorphisme de classes associ\'e \`a la multiplication par $n$ dans $NV$. La proposition \ref{produit} et le th\'eor\`eme \ref{nultor} permettent alors d'en d\'eduire que $\ker(\psi_n^{NV})=g^{-1}(\ker(\psi_n^A))$.

Nous obtenons d'autre part, en appliquant le lemme du serpent au diagramme pr\'ec\'edent, une suite exacte
$$
\begin{CD}
0 @>>> \mu_N @>>> V @>\phi>> NV @>>> 0.\\
\end{CD}
$$
De plus, les entiers $n$ et $N$ \'etant premiers entre eux, le morphisme $\phi$ induit un isomorphisme $\phi_0:V[n]\rightarrow NV[n]$. La proposition \ref{foncteurcor} affirme alors que
$$
\ker(\psi_n^V)=\phi^{-1}(\ker(\psi_n^{NV}))
$$
donc, en vertu de ce qui pr\'ec\`ede,
$$
\ker(\psi_n^V)=\phi^{-1}(g^{-1}(\ker(\psi_n^A)))=f^{-1}(\ker(\psi_n^A))
$$
ce qu'on voulait d\'emontrer.
\end{proof}


\vskip 2cm

Jean Gillibert
\smallskip

The University of Manchester

Alan Turing Building

Oxford Road

Manchester M13 9PL

United Kingdom

\bigskip

\texttt{jean.gillibert@manchester.ac.uk}


\begin{thebibliography}{12}

\bibitem[A1]{ref1} \textsc{A. Agboola}, {\it A geometric description of the class invariant homomorphism}, J. Th\'eor. Nombres Bordeaux {\bf 6} (1994), 273--280.
\bibitem[A2]{ref2} ---------, {\it Torsion points on elliptic curves and Galois module structure}, Invent. Math. {\bf 123} (1996), 105--122.
\bibitem[A3]{ref3} ---------, {\it On primitive and realisable classes}, Compositio Math. {\bf 126} (2001), 113--122.
\bibitem[B-M]{ref4} \textsc{H. Bass} and \textsc{M. P. Murthy}, {\it Grothendieck groups and Picard groups of abelian group rings}, Annals of Math. {\bf 86} (1967), 16--73.
\bibitem[BLR]{ref5} \textsc{S. Bosch}, \textsc{W. L\"utkebohmert} and \textsc{M. Raynaud}, {\it N\'eron Models}, Ergeb. Math. Grenzgeb. (3), vol. 21 (Springer, Berlin-Heidelberg-New York, 1990).
\bibitem[CN-T]{ref6} \textsc{P. Cassou-Nogu\`es} et \textsc{M. J. Taylor}, {\it Structures galoisiennes et courbes elliptiques}, J. Th\'eor. Nombres Bordeaux {\bf 7} (1995), 307--331.
\bibitem[Fr]{ref7} \textsc{A. Fr\"ohlich}, {\it Galois module structure of algebraic integers}, Ergeb. Math. Grenzgeb. (3), vol. 1 (Springer, Berlin-Heidelberg-New York, 1983).
\bibitem[G1]{ref8} \textsc{J. Gillibert}, {\it Invariants de classes : le cas semi-stable}, Compositio Math. {\bf 141} (2005), 887--901.
\bibitem[G2]{ref9} \textsc{J. Gillibert}, {\it Vari\'et\'es ab\'eliennes et invariants arithm\'etiques}, Annales de l'Institut
Fourier {\bf 56} (2006), 277--297.
\bibitem[G3]{ref10} \textsc{J. Gillibert}, {\it Invariants de classes : exemples de non-annulation en dimension sup\'erieure}, Mathematische Annalen {\bf 338} (2007), 475--495.
\bibitem[Gr]{ref11} \textsc{C. Greither}, {\it Cyclic Galois extensions of commutative rings}, Lecture Notes in Mathematics, vol. 1534 (Springer, Berlin-Heidelberg-New York, 1992).
\bibitem[EGA II]{ref12} \textsc{A. Grothendieck} et \textsc{J. Dieudonn\'e}, {\it \'El\'ements de g\'eom\'etrie alg\'ebrique, chapitre II : \'Etude globale \'el\'ementaire de quelques classes de morphismes}, Publ. Math. Inst. Hautes \'Etudes Sci. {\bf 8} (1961).
\bibitem[SGA 3]{ref13} \textsc{A. Grothendieck} et \textsc{M. Demazure}, {\it Sch\'emas en groupes}, Lecture Notes in Mathematics, vols. 151, 152, 153 (Springer, Berlin-Heidelberg-New York, 1970).
\bibitem[SGA 4]{ref14} \textsc{A. Grothendieck}, \textsc{M. Artin} et \textsc{J. L. Verdier}, {\it Th\'eorie des topos et cohomologie \'etale des sch\'emas}, Lecture Notes in Mathematics, vols. 269, 270 (Springer, Berlin-Heidelberg-New York, 1972).
\bibitem[SGA 7]{ref15} \textsc{A. Grothendieck}, {\it Groupes de monodromie en g\'eom\'etrie alg\'ebrique}, Lecture Notes in Mathematics, vol. 288 (Springer, Berlin-Heidelberg-New York, 1972).
\bibitem[Hil]{ref16} \textsc{D. Hilbert}, {\it Die Theorie der Algebraischen Zalhk\"orper}, Jahresbericht der Deutschen Mathematiker-Vereinigung {\bf 4} (1897), 175--546.
\bibitem[Mi80]{ref17} \textsc{J. S. Milne}, {\it \'Etale Cohomology}, Princeton Math. Ser., vol. 33 (Princeton University Press, 1980).
\bibitem[P]{ref18} \textsc{G. Pappas}, {\it On torsion line bundles and torsion points on abelian varieties}, Duke Math. J. {\bf 91} (1998), 215--224.
\bibitem[S-T]{ref19} \textsc{A. Srivastav} and \textsc{M. J. Taylor}, {\it Elliptic curves with complex multiplication and Galois module structure}, Invent. Math. {\bf 99} (1990), 165--184.
\bibitem[T81]{ref20} \textsc{M. J. Taylor}, {\it On Fr\"ohlich's conjecture for rings of integers of
tame extensions}, Invent. Math. {\bf 63} (1981), 41--79.
\bibitem[T88]{ref21} \textsc{M. J. Taylor}, {\it Mordell-Weil groups and the Galois module structure of rings of integers}, Illinois J. Math. {\bf 32} (1988), 428--452.
\bibitem[W]{ref22} \textsc{W. C. Waterhouse}, {\it Principal homogeneous spaces and group scheme extensions}, Trans. Am. Math. Soc. {\bf 153} (1971), 181--189.

\end{thebibliography}
\end{document}